\documentstyle[11pt,amsfonts]{article}
\def\pagebrE{\mbox{}}
\def\pagebrA{\pagebr}

\renewcommand{\epsilon}{\varepsilon}
\newtheorem{thm}{Theorem}[section]
\newtheorem{lem}[thm]{Lemma}
\newtheorem{prop}[thm]{Proposition}
\newtheorem{df}[thm]{Definition}
\newtheorem{cor}[thm]{Corollary}
\newtheorem{rem}[thm]{Remark}

\newenvironment{proof}{\medskip
\noindent {\bf Proof.}}{\hfill \rule{.5em}{1em}\mbox{}\bigskip}

\def\R{{\Bbb R}}
\def\C{{\Bbb C}}
\def\CP{{\Bbb C} {\Bbb P}}
\def\Q{{\Bbb Q}}

\def\F{{\Bbb F}}
\def\O{{\cal O}}
\def\FV{{\frak F}_V}
\def\so{{\frak so}}
\def\sl{{\frak sl}}
\def\g{{\frak g}}

\def\t{{\frak t}}
\def\L{{\frak L}}
\def\X{{\frak X}}
\def\gl{{\frak gl}}
\def\om{\omega}
\def\th{\theta}
\def\w{\wedge}
\def\und{\underline}
\def\lbr{<\!}
\def\rbr{\!>}
\def\l{\left}
\def\r{\right}
\def\ra{\rightarrow}
\def\a{{\bf a}}
\def\b{{\bf b}}
\def\c{{\bf c}}
\def\hook{\mbox{}\begin{picture}(10,10)\put(1,0){\line(1,0){7}}
  \put(8,0){\line(0,1){7}}\end{picture}\mbox{}}

\def\be{\begin{equation}}
\def\ee{\end{equation}}
\def\bi{\begin{itemize}}
\def\ei{\end{itemize}}
\def\ba{\begin{array}}
\def\ea{\end{array}}
\def\bea{\begin{eqnarray}}
\def\eea{\end{eqnarray}}
\def\linebr{\linebreak}
\def\pagebr{\pagebreak}

\title{\bf{On the Incompleteness of Berger's List of Holonomy
Representations}}
\author{Quo-Shin Chi\thanks{Supported in part by NSF Grant DMS
9301060}\\ Sergey A. Merkulov\\ Lorenz J. Schwachh\"{o}fer}
\date{August 30, 1995}
\begin{document}
\maketitle
\pagenumbering{arabic}

\begin{abstract}
In 1955, Berger \cite{Ber} gave a list of irreducible reductive
representations which can occur as the holonomy of a torsion-free
affine connection. This list was stated to be complete up to possibly a
finite number of missing entries.

In this paper, we show that there is, in fact, an infinite family of
representations which are missing from this list, thereby showing the
incompleteness of Berger's classification. Moreover, we develop a
method to construct torsion-free connections with prescribed holonomy,
and use it to give a complete description of the torsion-free affine
connections with these new holonomies. We also deduce some striking
facts about their global behaviour.
\end{abstract}

\section{Introduction}

An affine connection is one of the basic objects of interest in
differential geometry. It provides a simple and invariant way of
transferring information from one point of a connected manifold $M$ to
another and, not surprisingly, enjoys lots of applications in many
branches of mathematics, physics and mechanics. Among the most
informative characteristics of an affine connection is its (restricted)
holonomy group which is defined, up to conjugacy, as the subgroup of
$Gl(T_tM)$ consisting of all automorphisms of the tangent space $T_tM$
at $t \in M$ induced by parallel translations along $t$-based loops in
$M$.

Which reductive Lie groups $G$ can be irreducibly acting holonomies of
affine connections?

By a result of Hano and Ozeki \cite{HO}, {\em any} (closed) Lie group
representation $G \subseteq Gl(V)$ can be realized in this way. The
same question, if posed in the subclass of {\em torsion-free} affine
connections, has a very different answer. Long ago, Berger \cite{Ber}
presented a very restricted list of possible holonomies of a
torsion-free affine connection which, as he claimed, is complete up to
a finite number of missing terms. His list is separated into two
parts.  The first part corresponds to the holonomies of {\em metric}
connections, the second part to the {\em non-metric} ones. While Berger
gave detailed arguments for the proof of the metric part, the proof of
the second part was omitted.

The list of metric connections has been studied extensively in the
intervening years. In fact, it is by now well-known which entries of
this list actually {\em do} occur as holonomies of torsion-free
connections, and how the holonomy relates to both the geometry and the
topology of the underlying manifold. See \cite{Al}, \cite{Br1},
\cite{Br2}, \cite{J}, \cite{Si}, as well as the surveys in \cite{Sa}
and \cite{Bes}.

Despite the lack of proof, the second part of Berger's holonomy list
seems to have been generally accepted as correct. Even when Bryant
found examples of holonomies which are not on Berger's list \cite{Br3},
he called them {\em exotic holonomies}, suggesting that such holonomies
should exist in very special dimensions only, and therefore should be
thought of as analogous to the {\em exceptional holonomies} in the
metric case. Further examples of exotic holonomies were found in
\cite{CS}.

However, it is the subject of the present article to prove that, even
up to finitely many  missing terms, Berger's list is still incomplete.
This is done by proving the existence of an infinite family of
irreducible representations which are not on this list, yet do occur as
holonomy of torsion-free connections. These representations are:
\be \label{replist} \ba{llll}
Sl(2,\C) SO(n,\C), & {\mbox{acting on}} & {\R^{8n}} \cong \C^2
\otimes \C^n, & \mbox{where $n \geq 3$},\\
Sl(2,\R) SO(p,q), & {\mbox{acting on}} & {\R^{2(p+q)}} \cong \R^2
\otimes \R^{p+q}, & \mbox{where $p + q \geq 3$},\\
Sl(2,\R) SO(2,\R), & {\mbox{acting on}} & {\R^4} \cong \R^2
\otimes \R^2.\\
\ea \ee

\begin{thm} \label{thm:ConnectionsExist} All representations in {\em
(\ref{replist})} occur as holonomies of torsion-free connections which
are not locally symmetric.
\end{thm}

These candidates were first discovered by twistor theoretical methods
which shall be described in section \ref{sec:twistor}. In fact, the
formal existence of holomorphic torsion-free connections whose holonomy
is listed in the first family of (\ref{replist}) can be shown, in
principle, by deformation theory only.

We shall, however, present a different, technically simple method to
prove their existence in section \ref{sec:construct}. This new approach
will allow us to treat all representations in (\ref{replist})
simultaneously, and to assert some global geometric properties and
rigidity of connections with these holonomies which seem hard to
achieve otherwise. In fact, we shall classify {\em all} connections
with the above holonomies by showing that they all arise from this
construction. This implies, for example, the following:

\begin{thm} \label{thm:summary} Let $G \subseteq Gl(V)$ be one of the
representations in {\em (\ref{replist})} other than \linebr $Sl(2,\R)
SO(2,\R)$. Then the following hold.

\bi \item[{\em (1)}] Any torsion-free connection whose holonomy is
contained in $G$ is analytic.
\item[{\em (2)}] The moduli space of torsion-free connections whose
holonomy is contained in $G$ is finite dimensional. In particular, the
2nd derivative of the curvature at a single point completely determines
the connection.
\item[{\em (3)}] If $n \equiv 0,1 \bmod 4$ then any torsion-free
connection whose holonomy is contained in $G$ admits a non-trivial
infinitesimal symmetry, i.e. a vector field whose flow preserves the
connection.
\item[{\em (4)}] If $G \neq Sl(2,\R) SO(n,\R)$, then a torsion-free
connection whose holonomy is contained in $G$ is geodesically
incomplete, unless the connection is locally symmetric. If the latter
is the case, then the holonomy is a proper subgroup of $G$.
\ei
\end{thm}

The final step towards the proof Theorem \ref{thm:ConnectionsExist} was
accomplished when the second and third author met at the conference
{\em Geometry and Physics} in Aarhus, Denmark. We are grateful to the
organizers of this conference for inviting us; the third author wishes
to thank Professor Friedrich Hirzebruch for providing financial support
for his participation in that conference.

\section{Connections and Legendre moduli spaces} \label{sec:twistor}

Let $G \subseteq Gl(V)$ be an effective irreducible representation of a
complex connected reductive Lie group $G$ on the finite dimensional
complex vector space $V$. Without further comment, we shall always
regard $G$ as a Lie group with a fixed representation on $V$. Clearly,
$G$ also acts on $V^*$ via the dual representation. Let $\tilde{X}$ be
the $G$-orbit of a highest weight vector in $V^*$. Then the quotient
$X:= \tilde{X}/{\C}^*$ is a generalized flag variety which is
canonically embedded into ${\Bbb P}(V^*)$. In fact, $X=G_s/ P $, where
$G_s$ is the semi-simple part of $G$, and where $P$ is the parabolic
subgroup of $G_s$ leaving  the highest weight vector invariant up to a
scalar. Let $L_X$ be the restriction of the hyperplane bundle ${\cal
O}(1)$ on ${\Bbb P}(V^*)$ to the submanifold $X\hookrightarrow {\Bbb
P}(V^*)$. Then $V \cong H^0(X,L_X)$, and the Lie algebra $\g$ of $G$ is
contained in the Lie algebra
\[
H^0\l(X,L_X\otimes \l(J^1L_X\r)^*\r) \cong H^0(X,TX) \oplus {\C}.
\]
Therefore, to the representation $G \subseteq Gl(V)$, we canonically
associate a pair $(X,L_X)$, consisting of a generalized flag variety
$X$ and an ample line bundle $L_X$ on $X$. Most of the relevant
information about $G$ can be restored from $(X,L_X)$. The very few
cases when important information gets lost by the transition from
$(X,L_X)$ back to $G \subseteq Gl(V)$ are listed in \cite{Stei}.

A crucial step in deciding whether $G$ can occur as the holonomy of a
torsion-free affine connection is the computation of the formal
curvature space $K(\g)$, which is defined as the kernel of the
composition
\[
\Lambda^2 V^*\otimes {\g} \ra \Lambda^2 V^* \otimes V^*\otimes V \ra
\Lambda^3 V^* \otimes V.
\]
Namely, it follows from the {\em Ambrose-Singer Holonomy Theorem}
\cite{AS} that if there is a {\em proper} subalgebra $\g' \subset \g$
such that $K(\g) \subseteq \Lambda^2 V^* \otimes \g' \subset \Lambda^2
V^* \otimes \g$, then $G$ {\em cannot} be holonomy of a torsion-free
connection. This condition is called {\em Berger's first criterion},
and was used for the original classification in \cite{Ber}.

It is desirable to interpret $K(\g)$ in terms of the associated pair
$(X,L_X)$, because one will then be able to make use of the powerful
tools of complex analysis (such as, say, the vanishing theorems) to
attack the Berger classification problem. The twistor theory developped
in \cite{Me} does indeed provide us with such an interpretation. We
shall briefly describe this twistor construction, and use it to
calculate $K(\g)$ for the representations in (\ref{replist})
explicitly. For details and more general statements, we refer to
\cite{Me}.

Let $Y$ be a complex $(2n + 1)$-dimensional manifold. A {\em complex
contact structure}\, on $Y$ is a rank $2n$ holomorphic subbundle
$D\subseteq TY$ of the holomorphic tangent bundle to $Y$ such that the
Frobenius form
\begin{eqnarray*}
\Phi: D \times D & \longrightarrow & TY/D\\
(v,w) & \longmapsto & [v,w]\bmod D
\end{eqnarray*}
is non-degenerate. Define the contact line bundle $L$ by the exact
sequence
\[
0 \longrightarrow D \longrightarrow TY
\stackrel{\theta}{\longrightarrow} L \longrightarrow 0.
\]
One can easily verify that the maximal non-degeneracy of the
distribution $D$ is equivalent to the fact that the above defined
``twisted'' 1-form $\theta \in H^0(Y, L\otimes \Omega^1 M)$ satisfies
the condition
\[
\theta\wedge (d\theta)^n \neq 0.
\]

A complex $n$-dimensional submanifold $X$ of the complex contact
manifold $Y$ is called a {\em Legendre submanifold}\, if $TX \subseteq
D$.

Suppose $X \hookrightarrow Y$ is a compact complex Legendre submanifold
of a complex contact manifold $(Y, L)$, and denote the restriction
$L|_X$ by $L_X$. If $H^1(X, L_X) = 0$, then there exists a complete
family $\{X_t \hookrightarrow Y\ |\ t \in M \}$ of compact complex
Legendre submanifolds which is obtained from $X$ by all possible
holomorphic Legendre deformations of $X$ inside $Y$. The parameter
space $M$ is an $h^0(X, L_X)$-dimensional complex manifold, called the
{\em Legendre moduli space}. Moreover, if two more cohomology groups
vanish, namely if
\be \label{eq:cohomsvanish}
H^0\l(X, L_X \otimes S^2\l(J^1 L_X\r)^*\r) = H^1\l(X, L_X \otimes
S^2\l(J^1 L_X\r)^*\r) = 0,
\ee
then the Legendre moduli space $M$ comes equipped not only with an
induced complex manifold structure, but also with a uniquely induced
holomorphic torsion-free affine connection whose curvature tensor is a
field on $M$ with values in the cohomology group $H^1\l(X, L_X \otimes
S^3\l(J^1 L_X\r)^*\r)$. In fact, {\em any} holomorphic torsion-free
affine connection with reductive irreducibly acting holonomy group is
an induced connection on an appropriate Legendre moduli space
\cite{Me}.

This suggests the following strategy to look for new exotic
holonomies:  find a pair $(X, L_X)$ consisting of a generalized flag
variety $X = G/P$ and an ample line bundle $L_X \ra X$ such that the
cohomology groups in (\ref{eq:cohomsvanish}) vanish, while $H^1\l(X,
L_X \otimes S^3\l(J^1 L_X\r)^*\r)$ is non-zero. Then the twistor theory
guarantees that there is a natural injection
\be \label{inj}
\imath : H^1\l(X, L_X \otimes S^3\l(J^1 L_X\r)^*\r) \longrightarrow
\Lambda^2 V^* \otimes {\frak g},
\ee
where ${\g}:= H^0\l(X,L_X\otimes \l(J^1 L_X\r)^*\r)$ and $V:=
H^0(X,L_X)$, whose image equals $K(\g)$. In particular, $K(\g) \cong
H^1\l(X, L_X \otimes S^3\l(J^1 L_X\r)^*\r)$ as a $\g$-vector space. We
let $K_0(\g) \subseteq K(\g)$ be the set of elements with {\em full
curvature}, i.e.
\be \label{eq:fullcurvature}
K_0(\g) := \l\{ \rho \in K(\g)\ \l|\ \l< \l\{ \rho(x,y)\ |\ x,y \in V
\r\} \r> = \g \r. \r\}.
\ee
One can show that either $K_0(\g) = \emptyset$ or $K_0(\g)$ is {\em
dense in} $K(\g)$.

As a particular example, let us consider the reducible homogeneous
manifold $X =$ \linebr $\CP_1 \times \Q_m$, where $\Q_m
\hookrightarrow \CP_{m+1}$ is the non-degenerate quadric. We define the
line bundle
\[
L_X := \pi_1^*(\O_{\CP_1}(1)) \otimes \pi_2^*( \O_{
\CP_{m+1}}(1)|_{\Q_m}),
\]
where $\pi_1: X \ra \CP_1$ and $\pi_2: X \ra \Q_m$ are the canonical
projections. Then
\[
V:= H^0(X,L_X) = W_2\otimes W_n,
\]
where $W_2 := H^0(\CP_1, \O_{\CP_1}(1))$ and $W_n := H^0(\Q_m,
\O_{\CP_{m+1}}(1)|_{\Q_m})$ are vector spaces of dimensions $2$ and $n
:= m+2$, respectively, and
\[
\g:= H^0\l(X,L_X\otimes \l(J^1 L_X\r)^*\r) = {\g}_0 \oplus {\C},
\]
where ${\g}_0 = \sl(W_2) \oplus \so(W_n)$. Note that $W_2$ and $W_n$
carry a $\g_0$-invariant area form $\lbr\ ,\ \rbr$ and a
$\g_0$-invariant inner product $(\ ,\ )$, respectively. It is  easy to
check that for $m \geq 2$, i.e. $n \geq 4$, the cohomology groups in
(\ref{eq:cohomsvanish}) vanish, while
\be \label{eq:Kg-is-H1}
K(\g) \cong H^1\l(X, L_X \otimes S^3\l(J^1L_X\r)^*\r) \cong S^2(W_2)
\oplus \Lambda^2 W_n.
\ee
A calculation shows that $\Lambda^2 V^* \otimes \g$ contains only two
summands isomorphic to $S^2(W_2)$ and three isomorphic to $\Lambda^2
W_n$. After that, it is easy to obtain explicit expressions for the
elements of $K(\g)$. As it turns out, all elements of $K(\g)$ are
contained in $\Lambda^2 \otimes \g_0 \subset \Lambda^2 \otimes \g$. It
follows that {\em every torsion-free connection whose holonomy is
contained in $Gl(2,\C) SO(n,\C)$ for $n \geq 3$ must be contained in
$Sl(2,\C) SO(n,\C)$}.

In order to present the description of $K(\g)$, we use the
identifications $\sl(W_2) \cong S^2(W_2)$ and $\so(W_n) \cong \Lambda^2
W_n$, given by the identities
\be \label{eq:LieIsos} \ba{cccccc} (e_1 e_2)
\cdot e_3 & := & \lbr e_1, e_3 \rbr e_2 & + & \lbr e_2, e_3 \rbr e_1, &
\mbox{and}\\
(x_1 \w x_2) \cdot x_3 & := & (x_1, x_3) x_2 & - & (x_2, x_3) x_1
\ea \ee
for all $e_i \in W_2$, $x_i \in W_n$.

\begin{prop} \label{prop:compute-Kg} Let $G \subseteq Gl(V)$ be one of
the representations in (\ref{replist}) other than \linebr $Sl(2,\R)
SO(2,\R)$, and let $\g \subseteq \gl(V)$ be its Lie algebra. For $A \in
\sl(2,\C)$ and $M \in \so(n,\C)$ ($A \in \sl(2,\R)$ and $M \in
\so(p,q)$, respectively), define $\rho_{A+M}: \Lambda^2 V \ra \g$ by
\be \label{eq:curvature} \ba{ll}
\rho_{A + M} (e_1 \otimes x_1, e_2 \otimes x_2) := & \lbr e_1, e_2 \rbr
((x_1,x_2) (A + M)\ +\ (x_1 \w M x_2 + x_2 \w M x_1))\\
& +\ (x_1, M x_2) e_1 e_2\ -\ \lbr Ae_1, e_2 \rbr x_1 \w x_2.
\ea \ee
Then $\rho_{A+M} \in K(\g)$, and the map
\[
\ba{lcll} \rho: & \g & \longrightarrow & K(\g)\\
& A + M & \longmapsto & \rho_{A+M} \ea
\]
is a $G$-equivariant isomorphism.
\end{prop}

\begin{proof} It is straightforward to verify that $\rho_{A+M} \in
K(\g)$ and that $\rho$ is a $G$-equivariant injection. If $G = Sl(2,\C)
SO(n,\C)$ with $n \geq 4$, then the surjectivity of $\rho$ follows from
(\ref{eq:Kg-is-H1}). This implies the surjectivity of $\rho$ for $G =
Sl(2,\R) SO(p,q)$ with $p+q \geq 4$ as well.

If $n=3$ or $p+q=3$, then the surjectivity follows from a direct
calculation.
\end{proof}

\begin{rem} \label{remark} {\em \bi
\item[(1)] If $G = Sl(2,\R) SO(2,\R)$ then $\rho$ is injective but not
surjective. In fact, $dim(K(\g)) = 9 > dim(\g)$ in this case.

\item[(2)] The case $G = Sl(2,\F) SO(3,\F)$ with $\F = \R$ or $\C$ was
treated in \cite{CS}, where this representation was called $H_{12}$.

Note that the twistor approach from above does not work for $n=3$;
indeed, the cohomology $H^1\l(X, L_X \otimes S^2\l(J^1 L_X\r)^*\r) \neq
0$, thus (\ref{eq:cohomsvanish}) is {\em not} satisfied. Nevertheless,
it follows from the results in \cite{CS} that every Legendre moduli
space with \linebr $X \cong \CP_1 \times \Q_1 \cong \CP_1 \times \CP_1$
and $L_X$ as above, i.e. $L_X \cong \O(1,2)$, {\em does} carry an
induced torsion-free connection.  This shows, in particular, that
(\ref{eq:cohomsvanish}) is not a necessary condition for the existence
of torsion-free induced connections on a Legendre moduli space (cf.
\cite{Me}).

\item[(3)] If $A,M \neq 0$ then $\rho_{A+M}$ is surjective. Thus, for
the representations in (\ref{replist}), $K_0(\g) \subseteq K(\g)$ is
open dense.
\ei} \end{rem}

By (3), it follows that if there exists a holomorphic torsion-free
affine connection $\nabla$ with ``generic'' curvature tensor, then its
holonomy will be the full group $Sl(2,\C) SO(n,\C)$.

If $(Y,L)$ is a complex contact manifold containing $X = \CP_1 \times
\Q_m$ as a Legendre submanifold, then the latter is stable under
arbitrary holomorphic deformations of the contact structure on $Y$.
Therefore, one possible way to prove the existence of the required
holomorphic affine connections is to use the Kodaira deformation theory
\cite{Ko} to construct a sufficiently ``generic'' contact manifold
$Y$.

We shall, however, use a different method for the existence proof,
which will be presented in the following section.

\section{Construction of torsion-free connections}
\label{sec:construct}

Let us briefly recall the definition and basic properties of a Poisson
manifold. For a more detailed exposition, see e.g. \cite{L} or
\cite{V}.

\begin{df} \label{def:Poisson} A Poisson structure on a differentiable
manifold $P$ is a bilinear map, called the Poisson bracket
\[
\{\ ,\ \}: \otimes^2 C^\infty(P, \R) \longrightarrow C^\infty(P, \R),
\]
satisfying the following identities:

\bi
\item[{\em (i)}] the bracket is skew-symmetric:
\[ \{f,g\} = -\{g,f\}, \]
\item[{\em (ii)}] the bracket satisfies the Jacobi identity:
\[ \{f,\{g,h\}\} + \{g,\{h,f\}\} + \{h,\{f,g\}\} = 0, \]
\item[{\em (iii)}] the bracket is a derivation in each of its
arguments:
\[ \{f,gh\} = \{f,g\} h + g \{f,h\}. \]
\ei
\end{df}

It is well-known that on every Poisson manifold $(P, \{\ ,\ \})$, there
exists a unique smooth bivector field $\Lambda \in \Gamma(P, \Lambda^2
TP)$ such that the Poisson bracket is given by
\be
\label{eq:bivector} \{f,g\} = \Lambda(df,dg).
\ee
We define the homomorphism $\Lambda^\#: T^*P \ra TP$ by the equation
\be \label{eq:defLambda} \ba{ll}
(\Lambda^\# df)(g) = \{f, g\} & \mbox{for all $f,g \in C^\infty(P,
\R)$.}
\ea \ee
The {\em half-rank at $p \in P$} of the Poisson structure is the
smallest integer $r$ such that
\[
\Lambda^{r+1}(p) = 0,
\]
and the {\em rank at $p \in P$} is twice the half-rank. It follows that
the rank at $P$ equals the rank of $\Lambda^\#_p: T_p^*P \ra T_pP$. The
Poisson structure is called {\em non-degenerate at $p$} if
$\Lambda^\#_p$ is an isomorphism, i.e. if the rank at $p$ equals the
dimension of $P$. In particular, if $P$ is non-degenerate at a point
then $P$ must be even dimensional, and the set of non-degenerate points
is open in $P$. If $P$ is non-degenerate {\em everywhere}, then there
is a natural symplectic 2-form $\sigma$ on $P$ such that $\Lambda^\#$
is precisely the index-raising map associated to $\sigma$. In fact, it
is well known that symplectic structures are in a natural one-to-one
correspondence with non-degenerate Poisson structures.

The {\em characteristic field} of the Poisson structure is the subset
of $TP$ given by
\[
{\cal C} = \Lambda^\#(T^*P).
\]
Thus, the dimension of ${\cal C}_p$ equals the rank at $p$. A {\em
characteristic leaf} $\ \Sigma \subseteq P$ is a submanifold for which
$T_p\Sigma = {\cal C}_p$ for all $p \in \Sigma$. From
(\ref{eq:defLambda}), it follows that the set of functions which vanish
on $\Sigma$ form a {\em Poisson ideal;} hence there is a naturally
induced Poisson structure on $\Sigma$. Clearly, this Poisson structure
on $\Sigma$ is non-degenerate. Thus it follows that {\em every
characteristic leaf of a Poisson manifold carries a natural symplectic
structure}.

\begin{df} Let $(P, \{\ ,\ \})$ be a Poisson manifold. A symplectic
realization of $P$ is a symplectic manifold $(S, \sigma)$ and a
submersion
\[
\pi: S \longrightarrow P
\]
which is compatible with the Poisson structures, i.e.
\be \label{eq:PoissonMorph} \ba{ll}
\{\pi^*(f), \pi^*(g)\}_S = \pi^*(\{f,g\}) & \mbox{for all $f,g \in
C^\infty(P,\R)$,}
\ea \ee
where the Poisson bracket $\{\ ,\ \}_S$ on $S$ is induced by the
symplectic structure.
\end{df}

\begin{prop} \label{prop:SymplReal} If $p_0 \in P$ has locally constant
rank, i.e. the rank is constant on an open neighborhood $U$ of $p_0$,
then there is a local symplectic realization at $p_0$, i.e. a
symplectic realization $\pi: S \longrightarrow U'$ with $p_0 \in U'
\subseteq U$.
\end{prop}

\begin{proof} By Darboux's Theorem, there exists a local coordinate
system $(p_i, q_i, t_\alpha)$, \linebr $i = 1, \ldots, r$, $\alpha = 1,
\ldots, s$, in a neighborhood $U$ of $p_0$ such that the Poisson
bracket is given by
\[
\{f,g\} = \frac {\partial f}{\partial p_i} \frac {\partial g}{\partial
q_i} - \frac {\partial f}{\partial q_i} \frac {\partial g}{\partial
p_i}.
\]
Let $S := U \times \R^s$ with coordinates $u_\alpha$ on $\R^s$. Define
the symplectic 2-form \linebr $\sigma := dp_i \w dq_i + dt_\alpha \w
du_\alpha$ on $S$. Then it is easily verified that the natural
projection $\pi: S \ra U$ is a symplectic realization of $U$.
\end{proof}

We now turn to the construction of torsion-free connections via Poisson
structures. First of all, let us set up some notation.

Let $V$ be a finite dimensional vector space, $\g \subseteq \gl(V)$ a
Lie sub-algebra, and let $G \subseteq Gl(V)$ be the corresponding
connected Lie group. $G$ acts canonically on $V$, and on $\g$ via the
adjoint representation. This induces $G$-actions on all tensor powers
and direct sums of $\g$ and $V$ which we will call the {\em canonical
action} of $G$ on these spaces.

Recall that the {\em curvature space} $K(\g)$ is defined by the exact
sequence
\[
0 \longrightarrow K(\g) \longrightarrow \Lambda^2 V^* \otimes \g
\longrightarrow \Lambda^3 V^* \otimes V,
\]
where the latter map is the composition of the inclusion and the
skew-symmetrization map $\Lambda^2 V^* \otimes \g \ra \Lambda^2 V^*
\otimes V^* \otimes V \ra \Lambda^3 V^* \otimes V$.  Likewise, we
define the {\em 2nd curvature space} $K^1(\g)$ by the exact sequence
\[
0 \longrightarrow K^1(\g) \longrightarrow V^* \otimes K(\g)
\longrightarrow \Lambda^3 V^* \otimes \g,
\]
where again, the latter map is given by the composition of an inclusion
and skew-symmetrization, namely $V^* \otimes K(\g) \ra V^* \otimes
\Lambda^2 V^* \otimes \g \ra \Lambda^3 V^* \otimes \g$.  In other
words, $K(\g)$ and $K^1(\g)$ consist of those linear maps $\Lambda^2 V
\ra \g$ and $V \ra K(\g)$ which satisfy the 1st and 2nd Bianchi
identity, respectively.

Let $W := \g \oplus V$. Denote elements of $\g$ and $V$ by $A, B,
\ldots$ and $x, y, \ldots$, respectively, and elements of $W$ by $w,
w', \ldots$. We may regard $W$ as the semi-direct product of Lie
algebras, i.e. we define a Lie algebra structure on $W$ by the equation
\[
[A+x, B+y] := [A,B] + A \cdot y - B \cdot x.
\]
It is well-known \cite{L} that this induces a natural Poisson structure
on the dual space $W^*$. Now, we wish to perturb this Poisson
structure. For this, we need the
\pagebrA

\begin{df} A $C^\infty$-map $\phi: \g^* \ra \Lambda^2 V^*$ is called
admissible if
\bi \item[{\em (i)}] $\phi$ is $G$-equivariant,
\item[{\em (ii)}] $d\phi(p) \in K(\g) \subseteq \Lambda^2 V^* \otimes
\g$ for all $p \in \g^*$.
\ei
\end{df}
In order for condition (ii) to make sense, we use the natural
identification $T^*_p\g^* \cong \g$. Now, the following important
observation is easily proven.

\begin{prop} \label{prop:bracket}
Let $V$, $\g \subseteq \gl(V)$, $W$ and $K(\g)$ as above, and let
$\phi: \g^* \ra \Lambda^2 V^*$ be an admissible map. Let $\Phi := \phi
\circ pr$, where $pr: W^* \ra \g^*$ is the natural projection. Then the
following bracket on $W^*$ is Poisson:
\be \label{eq:bracket}
\{f,g\}(p) := p([A+x, B+y]) + \Phi(p)(x,y).
\ee
Here, $df_p = A + x$ and $dg_p = B + y$ are the decompositions of
$df_p, dg_p \in T^*_pW^* \cong W$.
\end{prop}

Note that for $\phi = 0$, we simply obtain the Poisson structure
induced by the Lie algebra structure on $W$.

Let us now consider a Poisson structure on $W^*$ induced by an
admissible map \linebr $\phi: \g^* \ra \Lambda^2 V^*$. Let $\pi: S \ra
U$ be a symplectic realization of the open subset $U \subseteq W^*$.
For each $w \in W$, we define the vector fields
\[
\eta_w := \Lambda^\#(w) \in \X(W^*),
\]
and
\[
\xi_w := \#(\pi^*(w)) \in \X(S),
\]
where $w \in W \cong T^*W^*$ is regarded as a 1-form on $W^*$. Since
$\pi$ preserves the Poisson structure (\ref{eq:PoissonMorph}), we have
\be \label{eq:pixi=eta} \ba{ll}
\pi_*(\xi_w) = \eta_w & \mbox{for all $w \in W$}.
\ea \ee
In contrast to the map $w \mapsto \eta_w$, the map $w \mapsto \xi_w$ is
{\em pointwise injective}. Thus, \linebr $\xi := \{\xi_w\ |\ w \in W \}
\subseteq TS$ is a distribution on $S$ whose rank equals the dimension
of $W$. For the bracket relations, we compute
\be \label{eq:xibracket} \ba{llll}
[\xi_A, \xi_B] & = & \xi_{[A,B]}\\
{[\xi_A, \xi_x]} & = & \xi_{A \cdot x}\\
{[\xi_x, \xi_y]}(s) & = & \xi_{d\Phi(p)(x,y)} & \mbox{ where $p =
\pi(s)$}.
\ea \ee

This implies, of course, that the distribution $\xi$ on $S$ is {\em
integrable}.  Moreover, the first equation in (\ref{eq:xibracket})
implies that the flow along the vector fields $\xi_A$ induces a local
$G$-action on $S$. Let $F \subseteq S$ be a maximal integral leaf of
$\xi$. Clearly, $F$ is $G$-invariant, and we can define a $W$-valued
coframe $\th + \om$ on $F$, where $\th$ and $\om$ take values in $V$
and $\g$, respectively, by the equation
\[
v_s = \xi_{(\om + \th) (v_s)} (s), \mbox{ all $v_s \in T_sF$}.
\]
The equations dual to (\ref{eq:xibracket}) then read
\be \label{eq:structureeqn} \ba{lll}
d\th & = & - \om \w \th\\
d\om & = & - \om \w \om - \pi^*(d\Phi) \circ (\th \w \th).
\ea \ee
Here, $d\Phi$ is regarded as a map with values in $K(\g) \subseteq
\Lambda^2 V^* \otimes \g$.

After shrinking $S$ and $U$ if necessary, we may assume that $M := F/G$
is a {\em manifold}. From (\ref{eq:structureeqn}) it follows that there
is a unique torsion-free connection on $M$ and a unique immersion
$\imath: F \hookrightarrow \FV$ into the $V$-valued coframe bundle
$\FV$ of $M$ such that $\th = \imath^*(\und\th)$ and $\om =
\imath^*(\und\om)$, where $\und\th$ and $\und\om$ are the tautological
and the connection 1-form on $\FV$, respectively. Clearly, the holonomy
of this connection is contained in $G$.

\begin{df} Let $\phi: \g^* \ra \Lambda^2 V^*$ be an admissible map.
Then a torsion-free connection which is obtained from the above
construction is called a Poisson connection induced by $\phi$.
\end{df}
We then get the following result.

\begin{thm} \label{th:ConnExist} Let $V$, $\g \subseteq \gl(V)$ and
$K(\g)$ be as before, and let $K_0(\g) \subseteq K(\g)$ be as in {\em
(\ref{eq:fullcurvature})}. Consider an admissible map $\phi: \g^* \ra
\Lambda^2 V^*$. Furthermore, suppose that
\[
\g^* \supseteq U_0 := (d\phi)^{-1} (K_0(\g)) \neq \emptyset.
\]
Then there exist Poisson connections induced by $\phi$ whose holonomy
representations are equivalent to $\g$. Moreover, if $\phi|_{U_0}$ is
not affine, then not all of these connections are locally symmetric.
\end{thm}

\begin{proof} Let $U^{reg} \subseteq U_0 \oplus V^* \subseteq W^*$ be
the subset of points for which the rank is locally constant. By upper
semi-continuity of the rank, $U^{reg}$ is open dense in $U_0 \oplus
V^*$.

Now Proposition \ref{prop:SymplReal} implies that there are symplectic
realizations $\pi: S \to U$, with open $U \subseteq U^{reg} \subseteq
W^*$. Then the above construction produces Poisson connections induced
by $\phi$ on some manifold $M = F/G$. By (\ref{eq:fullcurvature}),
(\ref{eq:structureeqn}) and the {\em Ambrose-Singer Holonomy Theorem}
\cite{AS}, the holonomy of this connection equals $\g$.

To show the last part, let us assume that {\em all} connections which
arise in this way are locally symmetric. Let $w := (p, q) \in U^{reg}$
with $p \in U_0$, $q \in V^*$. Then we may choose the symplectic
realization $\pi: S \ra U$ and $F \subseteq U$ such that $w \in
\pi(F)$. It is then easy to show by (\ref{eq:structureeqn}) that the
corresponding connection on $M := F/G$ is locally symmetric iff
$\L_{\xi_x}(\pi^*(d\Phi)) = 0$ for all $x \in V$. By
(\ref{eq:pixi=eta}) and because $\pi$ is a submersion, this is
equivalent to $\L_{\eta_x}(d\Phi) = 0 \mbox{ for all $x \in V$},$ or
$\L_{pr_*(\eta_x)}(d\phi) = 0 \mbox{ for all $x \in V$.}$ But now an
easy calculation shows that for all $A \in \g$,
\[
(pr_*(\eta_x)_w)(A) = -q(A \cdot x) = -\jmath(q \otimes x)(A),
\]
where $\jmath: V^* \otimes V \ra \g^*$ is the natural projection. Thus,
by our assumption, it follows that
\[
\L_{\jmath(q \otimes x)} (d\phi)_p = 0 \mbox{ for all $x \in V$, $(p,q)
\in U^{reg}$.}
\]
By density of $U^{reg}$, this equation holds for {\em all} $p \in U_0$,
$q \in V^*$, and since $\jmath$ is surjective, it follows that
\[
\L_\alpha(d\phi)_p = 0 \mbox{ for all $\alpha \in \g^*$, $p \in U_0$},
\]
i.e. $d\phi|_{U_0}$ is constant, hence $\phi|_{U_0}$ is affine.
\end{proof}

By Theorem \ref{th:ConnExist} it will suffice to address the question
of {\em existence} of admissible maps $\phi$ in order to construct
connections with prescribed holonomy.

Define the $k$-th jet space of $\g$ by
\[
J_k(\g) := \l(S^k(\g) \otimes \Lambda^2 V^*\r) \cap \l( S^{k-1}(\g)
\otimes K(\g) \r),
\]
where both are regarded as subspaces of $S^{k-1}(\g) \otimes \g \otimes
\Lambda^2 V^*$. Suppose there is a $G$-invariant element $\phi_k \in
J_k(\g)$. If we regard $\phi_k$ as a polynomial map of degree $k$,
$\phi_k: g^* \ra \Lambda^2 V^*$, then it follows that $\phi_k$ is
admissible. Conversely, given an {\em analytic} map $\phi: \g^* \ra
\Lambda^2 V^*$ with analytic expansion at $0 \in \g^*$
\[
\phi = \phi_0 + \phi_1 + \cdots,
\]
then it is straightforward to show that $\phi$ is admissible iff all
$\phi_k$ are, iff $\phi_k \in J_k(\g)^G$.

Consider an element $\phi_2 \in J_2(\g)^G$. On the one hand, we may
regard $\phi_2$ as an element of $\g \otimes K(\g)$, on the other hand,
it is easy to verify that also $\phi_2 \in V \otimes K^1(\g) \subseteq
V \otimes V^* \otimes K(\g)$. Thus, by the natural contractions,
$\phi_2$ induces $G$-equivariant linear maps
\be \label{eq:inducedmaps} \ba{llll}
\phi_2': & \g^* & \longrightarrow & K(\g)\\
\phi_2'': & V^* & \longrightarrow & K^1(\g).
\ea \ee

We shall now demonstrate the existence of torsion-free connections with
prescribed holonomy.

\begin{thm} \label{thm:PoissonExist} Let $G \subseteq Gl(V)$ be one of
the following representations:
\bi \item[{\em (i)}] For $\F = \R$ or $\C$, $G = Sl(2,\F)$ acting on
the space $V_3 \subseteq \F[x,y]$ of homogeneous polynomials in $x$ and
$y$ of degree $3$,
\item[{\em (ii)}] any of the representations in {\em (\ref{replist})}.
\ei

Then there is a $G$-invariant 2-form $\sigma \in \Lambda^2 V^*$ which
is unique up to a scalar. Also, $J_2(\g)$ is one-dimensional and acted
on trivially by $G$. The generic Poisson connection induced by the
admissible map \be \label{eq:admissible} \phi = \phi_2 + c \sigma, \ee
with $0 \neq \phi_2 \in J_2(\g)$ and some constant $c$, has full
holonomy $G$ and is not locally symmetric.
\end{thm}

This shows, in particular, that the representations in (\ref{replist})
do occur as holonomy representations, and hence proves Theorem
\ref{thm:ConnectionsExist}. Also, the representation in (i) is
precisely the representation $H_3$ which has already been shown to
occur as an exotic holonomy in \cite{Br3}.

\begin{proof} We shall present the proof for the representations in
(ii) only. In this case, the $G$-invariant symplectic form is
\[
\sigma(e_1 \otimes x_1, e_2 \otimes x_2) = \lbr e_1, e_2 \rbr
(x_1,x_2),
\]
where, as before, $\lbr \ ,\ \rbr$ and $(\ ,\ )$ denote the area form
on $\F^2$ and the inner product on $\F^n$, respectively, with $\F = \R$
or $\C$. We shall identify $\g$ and $\g^*$ via the Killing form $B$ on
$\g$, given by
\[
B(A + M, e_1 e_2 + x_1 \w x_2) = \lbr Ae_1, e_2 \rbr + (M x_1, x_2)
\mbox{ with $A \in \sl(\F^2)$ and $M \in \so(\F^n)$.}
\]
Here, once again we used the identifications $\sl(\F^n) \cong
S^2(\F^n)$ and $\so(\F^n) \cong \Lambda^2 \F^n$ from
(\ref{eq:LieIsos}). Now from the explicit description of $K(\g)$ in
Proposition \ref{prop:compute-Kg} it is straightforward to show that
$J_2(\g) = \l(\g \otimes K(\g)\r) \cap \l(S^2(\g) \otimes \Lambda^2
V\r)$ is one-dimensional and spanned by the element $\phi_2$, given by
\be \label{eq:phi2} \ba{ll}
\phi_2(A_1 + M_1, A_2 + M_2, & e_1 \otimes x_1, e_2 \otimes x_2) := \\
& \sigma(e_1 \otimes x_1, e_2 \otimes x_2) (B(A_1 + M_1, A_2 + M_2)) \\
& - (B(A_1, e_1 e_2) B(M_2, x_1 \w x_2) + B(A_2, e_1 e_2) B(M_1, x_1 \w
x_2))\\
& + \lbr e_1, e_2 \rbr ((M_1 x_1, M_2 x_2) + (M_1 x_2, M_2 x_1)).
\ea \ee

The isomorphism $\rho$ from (\ref{eq:curvature}) then coincides with
the map $\phi_2': \g^* \ra K(\g)$ from (\ref{eq:inducedmaps}), again
after identifying $\g$ and $\g^*$ via the Killing form $B$.  Moreover,
by (3) in Remark \ref{remark}, $K_0(\g) \subseteq K(\g)$ is {\em open
dense}. Then Theorem \ref{th:ConnExist} completes the proof.
\end{proof}

We can show even more. Namely, surprisingly enough, the converse of
Theorem \ref{thm:PoissonExist} is true:

\begin{thm} \label{th:allarePoisson} Let $G \subseteq Gl(V)$ be one of
the representations in Theorem \ref{thm:PoissonExist} other than
$Sl(2,\R) SO(2,\R)$. Then all torsion-free affine connections whose
holonomy is contained in $G$ are Poisson connections induced by the
admissible maps {\em (\ref{eq:admissible})}.
\end{thm}

This will follow immediately from the next Theorem, together with the
explicit description of the spaces $J_2(\g)$ in each case.

\begin{thm} \label{thm:allareinduced} Let $G \subseteq Gl(V)$ and $\g
\subseteq \gl(V)$ be as before. Suppose the following conditions are
satisfied:
\bi \item[{\em (i)}] $G$ is connected, reductive and acts irreducibly
on $V$,
\item[{\em (ii)}] there is a $\phi_2 \in J_2(\g)^G$ such that the
corresponding $G$-equivariant maps $\phi_2'$ and $\phi_2''$ from {\em
(\ref{eq:inducedmaps})} are isomorphisms.
\ei
Then every torsion-free affine connection whose holonomy is contained
in $G$ is a Poisson connection induced by a polynomial map
\[
\phi = \phi_2 + \tau,
\]
with $\phi_2 \in J_2(\g)$ from above and a (possibly vanishing)
$G$-invariant 2-form $\tau \in \Lambda^2 V^*$.
\end{thm}

For the proof, we shall need the following Lemma whose proof will be
given in the appendix.

\begin{lem} \label{lem:reps} Let $G \subseteq Gl(V)$ be an irreducible
representation of a connected, reductive Lie group $G$, and let $\g
\subseteq \gl(V)$ be the corresponding Lie algebra. If $\tau \in V^*
\otimes V^*$ satisfies the condition
\[
\tau(x, A \cdot y) = \tau(y, A \cdot x) \mbox{\ \ for all $x,y \in V$
and $A \in \g$,}
\]
then $\tau$ is skew-symmetric and hence a $G$-invariant 2-form.
\end{lem}

\medskip \noindent {\bf Proof of Theorem.} Let $F \subseteq \FV$ be a
$G$-structure on the manifold $M$ where $\FV \ra M$ is the $V$-valued
coframe bundle of $M$, and denote the tautological $V$-valued 1-form on
$F$ by $\th$. Suppose that $F$ is equipped with a torsion-free
connection, i.e. a $\g$-valued 1-form $\om$ on $F$. Since $\phi_2'$ is
an isomorphism, the {\em first and second structure equations} read
\be \label{eq:struct} \ba{ll}
d\th & = - \om \w \th\\
d\om & = - \om \w \om - 2 (\phi_2'(\a)) \circ (\th \w \th),
\ea \ee
where $\a: F \ra \g^*$ is a $G$-equivariant map. Differentiating
(\ref{eq:struct}) and using that $\phi_2''$ is an isomorphism yields
the {\em third structure equation} for the differential of $\a$:
\be \label{eq:struct3}
d\a = -\om \cdot \a + \jmath(\b \otimes \th),
\ee
for some $G$-equivariant map $\b: F \ra V^*$, where $\jmath:  V^*
\otimes V \ra \g^*$ is the natural projection. The multiplication in
the first term refers to the coadjoint action of $\g$ on $\g^*$. In
other words, (\ref{eq:struct3}) should be read as
\[ \ba{lll}
(\xi_A \a)(B) & = & \a([A,B])\\
(\xi_x \a)(B) & = & \b (B \cdot x).
\ea \]
Let us define the map $\c: F \ra V^* \otimes V^*$ by
\be \label{eq:struct4}
\c_p (x,y) := d\b(\xi_x) (y) - \phi_2(\a_p, \a_p, x, y).
\ee
Differentiation of (\ref{eq:struct3}) yields
\be \label{eq:cinvariance} \ba{ll}
\c_p(x, Ay) = \c_p(y, Ax) & \mbox{for all $x,y \in V$ and all $A \in
\g$.}
\ea \ee
Then by (i), (\ref{eq:cinvariance}) and Lemma \ref{lem:reps} we
conclude that $\c_p \in \Lambda^2 V^*$ is $G$-invariant. Moreover,
differentiation of (\ref{eq:struct4}) implies that $\xi_A(\c) = 0$ and
$(\xi_x \c)(y,z) = (\xi_y \c)(x,z)$ for all $A \in \g$ and $x,y,z \in
V$. Since $\c$ is skew-symmetric, it follows that
\[
d\c = 0,
\]
i.e. $\c_p \equiv \tau \in \Lambda^2 V^*$ is {\em constant}. Thus, the
$G$-equivariance of $\b$ and (\ref{eq:struct4}) yield
\be \label{eq:newstruct4}
d\b = -\om \cdot \b + \l( \a_p^2 \hook \phi_2 + \tau \r) \circ \th,
\ee
where $\hook$ refers to the contraction of $\a_p^2 \in S^2(\g^*)$ with
$\phi_2 \in S^2(\g) \otimes \Lambda^2 V^*$. In other words,
(\ref{eq:newstruct4}) should be read as
\[ \ba{lll}
(\xi_A \b)(y) & = & \b(A \cdot y)\\
(\xi_x \b)_p(y) & = & \phi_2(\a_p, \a_p, x, y) + \tau (x, y).
\ea \]

Let us now define the Poisson structure on $W^* = \g^* \oplus V^*$
induced by $\phi := \phi_2 + \tau$, and let $\pi := \a + \b: F \ra
W^*$. From (\ref{eq:struct3}) and (\ref{eq:newstruct4}) it follows that
$\pi_*(\xi_w) = \eta_w$ for all $w \in W$, and from there one can show
that, at least locally, the connection is indeed a Poisson connection
induced by $\phi$.
{\hfill \rule{.5em}{1em}\mbox{}\bigskip}

 From the complete characterization in Theorem \ref{th:allarePoisson},
we can deduce the following properties which summarize our discussion
so far:

\begin{cor} \label{cor:summary} Let $M$ be a manifold which carries a
torsion-free connection whose holonomy is contained in one of the
groups $G \subseteq Gl(V)$ from Theorem \ref{thm:PoissonExist} other
than $Sl(2,\R) SO(2,\R)$, and let $\phi = \phi_2 + c \sigma$ be the
admissible map which induces this connection. Then we have the
following.

\bi \item[{\em (1)}] The connection is analytic.
\item[{\em (2)}] The map $\pi:= \a + \b: F \ra W^*$ has constant even
rank $2k$ which we shall call the rank of the connection. $k=0$ iff the
connection is flat.
\item[{\em (3)}] $\pi(F)$ is contained in a $2k$-dimensional
characteristic leaf $\ \Sigma$ of the Poisson structure on $W^*$
induced by $\phi$. In particular, $\pi: F \ra \Sigma$ is a submersion
onto its image.
\item[{\em (4)}] Conversely, given a characteristic leaf $\Sigma
\subseteq W^*$ of maximal rank, then $\Sigma$ can be covered by open
neighborhoods $\{ U_\alpha \}$ such that there are Poisson connections
with $\pi(F_\alpha) = U_\alpha$.
\item[{\em (5)}] Let ${\frak s} \subseteq \X(F)$ be the Lie algebra of
infinitesimal symmetries of the connection, i.e. those vector fields
whose flows preserve the connection. Then $dim({\frak s}) = dim (W^*) -
2k$. In particular, if $n \equiv 0,1 \bmod 4$ then $dim({\frak s}) >
0$.
\item[{\em (6)}] The moduli space of torsion-free connections with any
of the above holonomies is finite dimensional. In particular, the 2nd
derivative of the curvature at a single point in $M$ completely
determines the connection on all of $M$.
\ei
\end{cor}

\begin{proof} (2) and (3) follow from the identity $\pi_*(\xi_w) =
\eta_w$. (1) and (4) are clear from the construction of the Poisson
connections and the analyticity of $\phi$, whereas (6) follows from the
structure equations in the proof of Theorem \ref{thm:allareinduced}.

To show (5), let $f: W^* \ra \F$ be a function which vanishes on
$\Sigma = \pi(F)$. Then it is easy to see that $\# \pi^*(df)$ is an
infinitesimal symmetry. It follows that $dim({\frak s}) \geq dim (W^*)
- 2k$. On the other hand, if $X \in {\frak s}$ then $\pi_*(X) = 0$,
hence $dim({\frak s}) \leq dim (W^*) - 2k$.
\end{proof}

Of course, (4) is not an optimal statement. One would like to show that
there are connections such that $\pi(F)$ is an {\em entire
characteristic leaf}. Also, our method does not prove the existence of
connections whose rank is not maximal. The difficulty is that, in
general, one cannot expect to have a {\em global} symplectic
realization $\pi: S \ra W^*$. In fact, even if we restrict to the
subset $W^{reg} \subseteq W^*$ where the Poisson structure has maximal
rank, then the obstruction for the existence of a global symplectic
realization is given by a class in $H^3_{rel}(W^{reg},{\cal F})$, where
$\cal F$ is the foliation by symplectic leafs \cite{V}.

Also, (5) is not an optimal statement either. Namely, as was shown in
\cite{CS}, if $n = 3$ then the dimension of the Lie algebra ${\frak s}$
of infinitesimal symmetries is at least 2. In fact, it seems likely
that $dim({\frak s}) > 0$ in {\em any} dimension, even though we do not
have a proof of this assertion.

For $H_3$-connections, all statements in Corollary \ref{cor:summary}
have been shown in \cite{Br3} by {\sc MAPLE} calculations. The same
kind of calculations was used in \cite{CS} to prove
Corollary~\ref{cor:summary} for $n = 3$. Our current approach, however,
seems more conceptual. In fact, it was a deeper understanding of the
structure of $H_3$-connections which led us to the construction of
Poisson connections presented in this paper.

There is another global result whose analogue for $H_3$-connections has
been demonstrated in \cite{Sch}:

\begin{thm} \label{thm:incomplete} Let $M$ be a manifold which carries
a torsion-free connection whose holonomy is contained in one of the
groups $G \subseteq Gl(V)$ from {\em (\ref{replist})} other than
$Sl(2,\R) SO(n,\R)$. Then $M$ is geodesically incomplete, unless the
connection is locally symmetric. If the latter is the case, then the
holonomy is a proper subgroup of $G$.
\end{thm}

\begin{proof} The holonomy of a locally symmetric connection must leave
the curvature tensor invariant. Since $K(\g)$ does not contain any
non-zero $G$-invariant element, the last statement follows.

Throughout this proof, we shall use the same notations and
identifications as in the proof of Theorem \ref{thm:PoissonExist}. In
particular, whenever it is convenient we shall regard $\a$ as a
$\g$-valued form, and write $\a_p = A_p + M_p$ with $A_p \in \sl(2,\F)$
and $M_p \in \so(n,\F)$.

$M$ is geodesically complete iff the vector fields $\xi_x \in \X(F)$
are complete for all $x \in V$, and we shall assume this from now on.
To prove the Theorem, we have to show that the connection is locally
symmetric, which is the case iff $\b \equiv 0$. If $M_p \equiv 0$ then
(\ref{eq:struct3}) implies that $\b \equiv 0$. Thus, we assume that at
some $p_0 \in F$, $M_{p_0} \neq 0$ and $\b_{p_0} \neq 0$ and shall
deduce a contradiction.

Let ${\cal N} := \{ y \in \F^n\ |\ (y,y) = 0 \}$. By the indefiniteness
of $(\ ,\ )$, ${\cal N}$ spans all of $\F^n$.  Therefore, $\b_{p_0}(e_1
\otimes y) \neq 0$ for some $e_1 \in \F^2, y \in {\cal N}$. Let $\nu :=
(M_{p_0} \cdot y, M_{p_0} \cdot y)$. If $\F = \R$ and $(\ ,\ )$ has
signature $(1, n-1)$, then it follows that $\nu \leq 0$. Therefore,
after possibly replacing $(\ ,\ )$ by its negative, we may assume that
$\nu \geq 0$ and, if $\F = \R$, the signature is $(p,q)$ with $p>1$. It
follows that there is a basis $\{x_1, \ldots, x_{n-2}, y, z\}$ of
$\F^n$ and $e_1, e_2 \in \F^2$ such that the following hold:

\be \label{eq:conditions} \ba{c}
(x_i, x_j) = \delta_i^j \epsilon_i \mbox{\ \ with\ \ }
\epsilon_i = \pm 1, \mbox{\ \ \ \ } \epsilon_1 = 1,\\
(x_i, y) = (x_i, z) = (y, y) = (z, z) = 0, \mbox{\ \ \ \ }
(y, z) = 1,\\
\lbr e_1, e_2 \rbr = 1, \mbox{\ \ \ \ }
(M_{p_0} \cdot y, M_{p_0} \cdot y) \geq 0, \mbox{\ \ and\ \ }
\b_{p_0}(e_1 \otimes y) \neq 0.
\ea \ee
Now, let us define
\[ \ba{lllll}
\xi_1 := \xi_{e_1 \otimes x_1}, &
\xi_2 := \xi_{e_2 \otimes y}, &
A_0 := -e_1^2 & \mbox{and} &
M_0 := 2 y \w x_1,
\ea \]
and the functions
\[ \ba{llll}
f_1(p) := \a_p(A_0), &
f_2(p) := \a_p(M_0) & \mbox{and} &
g(p) := - 2 \b_p(e_1 \otimes y).
\ea \]
 From (\ref{eq:phi2}), (\ref{eq:struct3}),
(\ref{eq:newstruct4}) and (\ref{eq:conditions}) we calculate
\be \label{eq:derivatives} \ba{c}
\xi_1(f_1) = \xi_2(f_2) = 0, \mbox{\ \ \ \ }
\xi_1(f_2) = \xi_2(f_1) = g,\\
\xi_1(g) = 2 f_1 f_2, \mbox{\ \ and\ \ }
\xi_2(g) = 4 (M_p \cdot y, M_p \cdot y).
\ea \ee
For $v \in \F^n$, we have $(v,v) = \sum_i \epsilon_i (v,x_i)^2 + (v,y)
(v,z)$. Thus,
\[
(M_p \cdot y, M_p \cdot y) = \sum_i \epsilon_i (M_p \cdot y, x_i)^2 =
\sum_i \epsilon_i (\a_p (y \w x_i))^2.
\]
We have $\a_p (y \w x_1) = B(M_p, y \w x_1) = \frac 12 f_2(p)$, and
from (\ref{eq:struct3}) and (\ref{eq:conditions}) we get that
$\xi_k(\a_p (y \w x_i)) = \b_p((y \w x_i) \cdot \xi_k) = 0$ for $i>1$.
Therefore,
\be \label{eq:squareplusconst} \ba{ll}
\xi_k \l( \xi_2(g) - f_2^2 \r) = 0 & \mbox{for $k = 1,2$.}
\ea \ee
For arbitrary constants $0 \neq c_1, c_2 \in \F$, define the vector
field $\xi$ and the function $f$ by
\[
\xi := c_1 \xi_1 + c_2 \xi_2 \mbox{\ \ and\ \ }
f := c_1^2 f_1 + 2 c_1 c_2 f_2.
\]
 From (\ref{eq:derivatives}) we compute that
\[ \ba{ll}
\xi^2(f) & = 3 c_1^2 c_2 \xi(g)\\
& = 3 c_1^2 c_2 (2 c_1 f_1 f_2 + c_2 \xi_2(g))\\
& = f^2 - c_1^2 \l( (c_1 f_1 - c_2 f_2)^2 - 3 c_2^2 (\xi_2(g) - f_2^2)
\r).
\ea \]
By (\ref{eq:squareplusconst}), $\xi \l( (c_1 f_1 - c_2 f_2)^2 - 3 c_2^2
(\xi_2(g) - f_2^2) \r) = 0$. It follows that along the flow line of
$\xi$, $f$ satisfies the differential equation
\be \label{eq:ODE}
y'' = y^2 + C,
\ee
where $C$ is a constant. By the assumption of completeness, $f$ must be
a {\em global} solution of (\ref{eq:ODE}). The following Lemma will be
shown in the appendix.
\pagebrA

\begin{lem} \label{lem:ODE} \bi
\item[{\em (i)}] The only holomorphic solutions of {\em (\ref{eq:ODE})}
which are defined on all of $\C$ are constants.
\item[{\em (ii)}] Suppose there is a real solution of {\em
(\ref{eq:ODE})} which is defined on all of $\R$. If for $x_0 \in \R$,
$y(x_0) > 0$ and $y'(x_0) \neq 0$, then $y''(x_0) < 0$.
\ei \end{lem}

We shall use this Lemma to get the desired contradiction. If $\F = \C$,
the Lemma implies that $3 c_1^2 c_2 g = \xi(f) \equiv 0$, contradicting
$g(p_0) = -2 \b_{p_0}(e_1 \otimes y) \neq 0$.

Consider now the case $\F = \R$. From (\ref{eq:derivatives}) and
(\ref{eq:squareplusconst}) we have $(\xi_2)^3(f_1) = 0$. Moreover,
$(\xi_2)^2(f_1)_{p_0} = 4 (M_{p_0} \cdot y, M_{p_0} \cdot y) \geq 0$
and $\xi_2(f_1)_{p_0} = g(p_0) \neq 0$. Therefore, since $\xi_2$ is
complete, there is a point $q_0$ on the flow line of $\xi_2$ through
$p_0$ which, in addition to (\ref{eq:conditions}), also satisfies
$f_1(q_0) > 0$.  W.l.o.g. we assume that $q_0 = p_0$.

Now $y(x_0) = f(p_0) = (c_1^2 f_1 + 2 c_1 c_2 f_2)(p_0)$ and $y''(x_0)
= \xi^2(f)_{p_0} = 3 c_1^2 c_2 (2 c_1 f_1 f_2 + c_2 \xi_2(g))(p_0)$.
Since $f_1(p_0) > 0$ and $\xi_2(g)_{p_0} \geq 0$ it follows that for a
suitable choice of $c_1, c_2$, both $y(x_0) > 0$ and $y''(x_0) \geq 0$
can be achieved. But $y'(x_0) = g(p_0) \neq 0$, hence again, we get a
contradiction from the Lemma.
\end{proof}

It is now clear that Theorem \ref{thm:summary} is simply a restatement
of parts of Corollary \ref{cor:summary} and Theorem
\ref{thm:incomplete}.
\pagebrE

\section{Appendix}

We shall now give the proof of two technical Lemmas.

\noindent {\bf Lemma \ref{lem:reps}} {\em Let $G \subseteq Gl(V)$ be an
irreducible representation of a connected, reductive Lie group $G$, and
let $\g \subseteq \gl(V)$ be the corresponding Lie algebra. If $\tau
\in V^* \otimes V^*$ satisfies the condition
\[
\tau(x, A \cdot y) = \tau(y, A \cdot x) \mbox{\ \ for all $x,y \in V$
and $A \in \g$,}
\]
then $\tau$ is skew-symmetric and hence a $G$-invariant 2-form.}

\begin{proof} Clearly, the Lemma is invariant under complexification,
hence we may assume that $G, \g$ and $V$ are complex. Also, if $\g_s
\subseteq \g$ is the semi-simple part of $\g$, then $\g_s$ acts
irreducibly on $V$ as well. Thus, we may assume w.l.o.g. that $\g$ is
semi-simple. Let
\[ \ba{lll}
\g = \t \oplus \bigoplus_{\alpha \in \Delta} \g_{\alpha} & \mbox{and} &
V = \bigoplus_{\lambda \in \Lambda} V_\lambda
\ea \]
be a Cartan and weight space decomposition of $\g$ and $V$,
respectively. We shall denote elements of $\t$, $\g_\alpha$ and
$V_\lambda$ by $A_0$, $A_\alpha$ and $x_\lambda$, respectively. Also,
we always write $\alpha, \beta, \ldots$ for roots, whereas $\lambda,
\mu, \ldots$ stand for weights.  Let $\lambda_{max} \in \t^*$ be the
maximal weight, and $x_{max} \in V_{\lambda_{max}}$ be a highest weight
vector.

The proof now proceeds as follows. Given $\tau$ as above, we shall
prove:
\[ \ba{ll}
\mbox{Step 1:} & \mbox{if $\g \cong \sl(2,\C)$ and $\lambda + \mu \neq
0$, then $\tau(x_\lambda, x_\mu) = 0$.}\\
\mbox{Step 2:} & \mbox{$\tau(x_\lambda, x_\mu) \mu = \tau(x_\mu,
x_\lambda) \lambda$ for all $\lambda, \mu \in \Lambda$. In particular,
$\tau(x_\lambda, x_\mu) = 0$ if $\lambda, \mu$}\\ & \mbox{are linearly
independent.}\\
\mbox{Step 3:} & \mbox{if $\lambda, \mu$ are scalar multiples of
$\lambda_{max}$ and $\lambda, \mu, \lambda + \mu \neq 0$ then}\\ &
\tau(x_\lambda, x_\mu) = \tau(x_\mu, x_\lambda) = 0.\\
\mbox{Step 4:} & \mbox{$\tau(x_{max}, x_0) = \tau(x_0, x_{max}) =
0$.}\\
\mbox{Step 5:} & \mbox{$\tau$ is skew-symmetric.}
\ea \]

{\em Proof of step 1:} Suppose $\g \cong \sl(2, \C)$, and $V \cong V_n$
is the (unique) irreducible $n + 1$-dimensional representation of $\g$.
Then it is an easy exercise to show that the only elements $\tau \in
V^* \otimes V^*$ which satisfy the above identity are $\tau = 0$ if $n$
is even, and a multiple of the symplectic $\g$-invariant 2-form on $V$
if $n$ is odd. From this, step 1 follows. The proof is omitted.

{\em Proof of step 2:} $\tau(x_\lambda, A_0 \cdot x_\mu) = \tau(x_\mu,
A_0 \cdot x_\lambda)$, hence $\mu(A_0) \tau(x_\lambda, x_\mu) =
\lambda(A_0) \tau(x_\mu, x_\lambda)$ for all $A_0 \in \t$. This proves
step 2.

{\em Proof of step 3:} Let $\lambda, \mu$ be as above, and let $\alpha
\in \Delta$ be linearly independent of $\lambda_{max}$. We compute
$\tau(x_\lambda, A_\alpha A_{-\alpha} \cdot x_\mu) = \tau(A_{-\alpha}
\cdot x_\mu, A_\alpha \cdot x_\lambda) = 0$ by step 2; indeed,
$A_{-\alpha} \cdot x_\mu \in V_{\mu - \alpha}$, $A_\alpha \cdot
x_\lambda \in V_{\lambda + \alpha}$, and $\mu - \alpha$, $\lambda +
\alpha$ are linearly independent if $\lambda + \mu \neq 0$. Likewise,
$\tau(x_\lambda, A_{-\alpha} A_\alpha \cdot x_\mu) = 0$, and hence $0 =
\tau(x_\lambda, [A_\alpha, A_{-\alpha}] \cdot x_\mu) = \mu ([A_\alpha,
A_{-\alpha}]) \tau(x_\lambda, x_\mu)$. Thus, we have $\lambda_{max}
([A_\alpha, A_{-\alpha}]) \tau(x_\lambda, x_\mu) = 0$ and similarly,
$\lambda_{max} ([A_\alpha, A_{-\alpha}]) \tau(x_\mu, x_\lambda) = 0$.
But if $\lambda_{max} ([A_\alpha, A_{-\alpha}]) = 0$ for all roots
$\alpha$ independent of $\lambda_{max}$, then it follows that $\g$ must
contain $\sl(2,\C)$ as a summand, and $\lambda_{max}$ is a multiple of
the root corresponding to this summand.

In this case, however, since the representation is faithful, it must be
equivalent to the action of $\sl(2,\C)$ on $V_n$. This observation
together with step 1 implies step 3.
\pagebrE

{\em Proof of step 4:} The equation $\tau(x_0, x_{max}) = 0$ follows
immediately from step 2.

Since the representation is irreducible, we have $V_0 = span
\{A_{-\alpha} \cdot V_\alpha\ |\ \alpha \in \Delta\}$. Thus, in order
to show step 4, we need to show
\be \label{eq:lemproof}
\tau(x_{max}, A_{-\alpha} \cdot x_\alpha) = 0 \mbox{ for all } \alpha
\in \Delta, x_\alpha \in V_\alpha.
\ee

If $\alpha \neq \lambda_{max}$ we have $\tau(x_{max}, A_{-\alpha} \cdot
x_\alpha) = \tau(x_\alpha, A_{-\alpha} \cdot x_{max}) = 0$ by steps 2
and 3; indeed, $A_{-\alpha} \cdot x_{max} \in V_{\lambda_{max} -
\alpha}$ and $\lambda_{max} - \alpha \notin \{ 0, -\alpha \}$. Thus,
(\ref{eq:lemproof}) follows in this case.

Next, suppose that $\alpha = \lambda_{max}$, and let $\beta$ and
$x_{\alpha + \beta}$ be such that $x_{max} = A_{-\beta} \cdot x_{\alpha
+ \beta}$. By irreducibility, this is always possible. Note that $\beta
\neq \alpha$, since by maximality of $\alpha = \lambda_{max}$, $2
\alpha$ is not a root. Then
\[ \ba{ll}
\tau(x_{max}, A_{-\alpha} \cdot x_\alpha) & =
\tau(x_\alpha, A_{-\alpha} A_{-\beta} \cdot x_{\alpha + \beta})\\ & =
\tau(x_\alpha, ([A_{-\alpha}, A_{-\beta}] +
A_{-\beta} A_{-\alpha}) \cdot x_{\alpha + \beta})\\ & =
\tau(x_{\alpha + \beta}, [A_{-\alpha}, A_{-\beta}] \cdot x_\alpha) +
\tau(A_{-\alpha} \cdot x_{\alpha + \beta}, A_{-\beta} \cdot x_\alpha)\\
& := I + II.
\ea \]
Now $x_{\alpha + \beta} \in V_{\alpha + \beta}$ and $[A_{-\alpha},
A_{-\beta}] \cdot x_\alpha \in V_{-\beta}$. Thus, by step 2 and the
fact that $\alpha \neq \beta$ it follows that $I = 0$. For $II$, we
have $A_{-\alpha} \cdot x_{\alpha + \beta} \in V_\beta$ and $A_{-\beta}
\cdot x_\alpha \in V_{\alpha - \beta}$. Now if $\beta$ is independent
of $\alpha$, step 2 implies that $II = 0$. If $\beta = -\alpha$ then
$A_{-\beta} \cdot x_\alpha \in V_{\alpha - \beta} = V_{2 \alpha} = 0$
by maximality of $\lambda_{max} = \alpha$, and again, $II = 0$
follows.  Thus, (\ref{eq:lemproof}) is shown in this case as well, and
step 4 follows.

{\em Proof of step 5:} Let $\tau^{sym}(x,y) := \tau(x,y) + \tau(y,x)$.
 From steps 2 through 4, it follows that $\lambda_{max} + \lambda \neq
0$ implies that $\tau(x_{max}, x_\lambda) = \tau(x_\lambda, x_{max}) =
0$. Also, if $x_{-max} \in V_{-\lambda_{max}}$, then by step 2,
$\tau^{sym}(x_{max}, x_{-max}) = 0$.

Thus, we conclude that $\tau^{sym}(x_{max},\und{\ \ }) = 0$. But this
is true for {\em any} choice of Cartan decomposition. Exploiting all
possible decompositions, we conclude
\[
\tau^{sym}(g \cdot x_{max}, \und{\ \ }) = 0 \mbox{ for all } g \in G.
\]
By irreducibility, the $G$-orbit of $x_{max}$ spans all of $V$, and
from there it follows that \linebr $\tau^{sym} = 0$, hence $\tau$ is
skew-symmetric.
\end{proof}

\noindent {\bf Lemma \ref{lem:ODE}} {\em Consider the differential
equation $y'' = y^2 + C$ {\em (\ref{eq:ODE})} for a constant $C$.
\bi \item[{\em (i)}] The only holomorphic solutions of {\em
(\ref{eq:ODE})} which are defined on all of $\C$ are constants.
\item[{\em (ii)}] Suppose there is a real solution of {\em
(\ref{eq:ODE})} which is defined on all of $\R$. If for $x_0 \in \R$,
$y(x_0) > 0$ and $y'(x_0) \neq 0$, then $y''(x_0) < 0$.
\ei}

\begin{proof} First of all, we integrate (\ref{eq:ODE}) to
\be \label{eq:ODE2}
(y')^2 = \frac 23 y^3 + 2 C y + C_1
\ee
for some constant $C_1$.

(i) Let $y$ be a holomorphic solution which is an entire function.
Define
\[ \ba{llll}
\phi: & \C & \longrightarrow & \CP_2\\ & x & \longmapsto &
[y(x) : y'(x) : 1].
\ea \]
By (\ref{eq:ODE2}), the image of $\phi$ is contained in the cubic $\cal
C$ given by $3 q^2 r - 2 p^3 - 6 C p r^2 - 3 C_1 r^3 = 0$, where
$p,q,r$ are the coordinates in $\CP_2$.

We calculate that $\cal C$ is a regular curve of genus 1 if $16 C^3 + 9
C_1^2 \neq 0$; it is a rational curve with a node if $16 C^3 + 9 C_1^2
= 0$ and $C \neq 0$, and it is a rational curve with a cusp if $C = C_1
= 0$.

Since $im(\phi)$ does not contain the point $[0:1:0] \in \cal C$, it
follows that in the first case $\phi$ maps $\C$ non-surjectively to a
genus 1 curve. By Picard's Theorem, this implies that $\phi$, and hence
$y$, is constant.

In the second case, we find a parametrization of $\cal C$ and compute
that there must be a function $t = t(x)$ such that
\[ \ba{lll}
y = \frac 3 {2C} \l( C_1 - 4 C^2 t^2 \r) & \mbox{and} &
y' = \frac{2 C} {C_1} \l( 3 C_1 - 8 C^2 t^2 \r) t.
\ea \]
But this implies that either $t \equiv 0$ or $t' = \frac{4 C^2}{3 C_1}
t^2 - \frac 12$.  It is straightforward to verify that the only global
solutions of the latter equation are constants.

In the third case $C = C_1 = 0$, a parametrization is given by $y = 6
t^2, y' = 12 t^3$ for some entire function $t$. But then $12 t t' = 12
t^3$, which again has no global non-constant solution.

(ii) Consider a global real solution $y$ and $x_0 \in \R$ as above.
After possibly replacing $y(x)$ by the solution $y(-x)$ and $x_0$ by
$-x_0$, we may assume that $y'(x_0) > 0$. If $y^2(x_0) + C = y''(x_0)
\geq 0$, then elementary calculus shows that $\lim_{x \ra \infty} y =
\infty$.  Thus, from (\ref{eq:ODE2}), we have for sufficiently large
$x$ that $(y')^2 > \frac 13 y^3$, and hence
\[
\l( y^{-\frac 12} \r)' = - \frac 12 y^{-\frac 32} y' <
- \frac 1 {2 \sqrt 3},
\]
contradicting $y^{-\frac 12} > 0$ for large $x$. Thus, $y''(x_0) \geq
0$ is impossible.
\end{proof}

\[ \small \ba{ll}
\mbox{\sc Quo-Shin Chi} &
\mbox{\sc Department of Mathematics, Campus Box 1146,}\\ &
\mbox{\sc Washington University, St. Louis, Mo 63110, USA}\\ &
\mbox{\rm chi@artsci.wustl.edu} \\

\mbox{\sc Sergey A.  Merkulov} &
\mbox{\sc Department of Pure Mathematics, Glasgow Univer-}\\ &
\mbox{\sc sity, 15 University Gardens, Glasgow G12 8QW, UK}\\ &
\mbox{\rm sm@maths.glasgow.ac.uk} \\

\mbox{\sc Lorenz J.  Schwachh\"{o}fer} &
\mbox{\sc Max-Planck-Institut f{\"u}r Mathematik,}\\ &
\mbox{\sc Gottfried-Claren-Stra{\ss}e 26, 53225 Bonn, Germany}\\ &
\mbox{\rm lorenz@mpim-bonn.mpg.de}
\ea
\]

\end{document}